\documentclass[10pt]{amsart}

\usepackage{amsxtra}
\usepackage{amssymb}

\newcommand\la{\langle}
\newcommand\ra{\rangle}

\newcommand\dd{{\mathfrak d}}
\renewcommand\aa{{\mathfrak a}}
\newcommand\qq{{\mathfrak q}}

\newcommand\rr{{\mathfrak r}}
\newcommand\hh{{\mathfrak h}} 
\newcommand\nn{{\mathfrak n}} 
\newcommand\ggo{{\mathfrak g}} 

\newcommand\ee{{\mathfrak e}}
\newcommand\aff{{\mathfrak{aff}}}

\newcommand\CC{\mathbb C}

\newcommand\RR{\mathbb R}

\newcommand\GL{{\mathrm{GL}}}

\newcommand\Aut{\operatorname{Aut}}

\newcommand\id{\operatorname{Id}}
\newcommand\SL{\operatorname{SL}}
\newcommand\Rg{\operatorname{im}}

\renewcommand\Re{\operatorname{Re}}
\renewcommand\Im{\operatorname{Im}}

\theoremstyle{plain}

\newtheorem{thm}{Theorem}[section]

\newtheorem{prop}[thm]{Proposition}
\newtheorem{cor}[thm]{Corollary}

\theoremstyle{definition}
\newtheorem{defn}[thm]{Definition}
\newtheorem{remark}[thm]{Remark}

\begin{document}
\title[complex, symplectic and K\"ahler structures on Lie groups]{complex, symplectic and K\"ahler structures \\ on four dimensional Lie groups}

\author{Gabriela Ovando}
\address{CIEM - Facultad de Matem\'atica, Astronom\'\i a y F\'\i sica,
Universidad Nacional de C\'or\-do\-ba, C\'or\-do\-ba~5000, Argentina}
\email{ovando@mate.uncor.edu}


\begin{abstract}In this work we deal with left invariant complex and symplectic structures  on simply connected four dimensional solvable real Lie groups. We search the general  form of such structures, when they exist and we make use of this information to determine all left invariant  K\"ahler structures. Finally, as an appendix we compute explicitly the cohomology over $\RR$ of the corresponding Lie algebras.
\end{abstract}

\thanks{{\it (2000) Mathematics Subject Classification}: 53C15, 53C55, 53D05, 22E25, 17B56 }

\maketitle

\section{Introduction}

The study of complex and symplectic manifolds has attracted the attention of many authors interested in different fields in mathematics and physics since complex and symplectic structures have proved to be an important tool in the description and geometrization of several phenomena. A complex manifold posseses a complex structure  on the underlying real manifold $M$, that is a differential tensor $J$ defined on $TM$ such that: (1) $J^2=-\id$ and (2) $N_J(X,Y) =0$ for all $X,Y \in \Xi(M)$, the so-called Nijenhuis condition. It is also known that the existence of such a tensor $J$ on a real manifold $M$, which satisfies condition (1) and (2) allows to construct complex coordinates on $M$.

If the manifold is homogeneous, the list of articles related with complex and/or symplectic questions is quite large. Existence conditions of symplectic forms on homogeneous spaces were found by Chu \cite{Ch} in the beginning 1970s. On the other hand homogeneous complex manifolds for example were classified in different stages \cite{O-R} \cite{Wi}.

From the point of view of Lie groups the existence problem of left invariant complex structures was treated by Samelson \cite{Sm} and Wang \cite{Wa} in the compact case and by Morimoto \cite{Mo} in the non compact case. Detailed classifications of left invariant complex structures when the Lie group is $\GL(2,\RR),\, U(2),\, \SL(3,\RR)$ were given by Sasaki \cite{Ss1}, \cite{Ss2}, and when $G$ is a simply connected solvable Lie algebra of dimension four by Snow \cite{Sn} and the author \cite{O}. The last two works were oriented to the study of the existence of left invariant complex structures and the determination of  the underlying complex manifold. The method used in these articles was based on the algebraic data related to the Lie algebra. In particular the techniques used for the explicite classifications in the last cases seem not to be adequate for the similar purpose in higher dimensions. Another point of view is treated in other kind of articles, which  deal with other features of complex and symplectic structures (just only in the solvable case see \cite{A-G-S} \cite{B-G} \cite{C-F-G-U}  \cite{F-G} \cite{G-R} \cite{Sl} or  for the study  of the moduli space of complex structures \cite{A-F-G-M}). 

The addition of an extra structure on the differential manifold such as (pseudo)-Riemannian metric or a symplectic structure compatible or not with  the complex structure introduces a new tool to be handled on the investigation of these objects. In this context the methods of rational homotopy theory have been applied successfully to symplectic geometry. These techniques were useful in attacking many geometric problems such as the construction of symplectic manifolds with no K\"ahler structure among others (see \cite{O-T} for a compilation on this topic and other articles such as \cite{A} \cite{A-G}  \cite{F-G-G} \cite{L} \cite{Mc} for example).

The purpose of this article is to study the space of left invariant complex and symplectic structures on solvable simply connected four dimensional real Lie groups (in fact, it is well known that semisimple Lie groups do not admit symplectic structures \cite{Ch}).
Based on a table presented in the next section, which  describes all solvable Lie algebras of dimension four (see \cite{A-B-D-O}) and by using explicit computations in the next two sections we describe the  general form of such a complex or symplectic structure. Note that particular cases, namely exact symplectic structures, were related with the existence of invariants \cite{Ca}.  
 All this information was used in the last section to search for K\"ahler structures. From our results we observe that not any case admits a complex and/or a symplectic structure. Moreover there are examples in which there exist a complex and a symplectic structure but not a K\"ahler structure. As an appendix we compute explicitely the real cohomology of the corresponding Lie algebras.

The author thanks I. Dotti for useful comments  and general supervision and specially M. Fern\'andez for her valuable contributions and for suggesting the subject of this article.

\

If $G$ is a simply connected Lie group then its Lie algebra will be denoted with greek letters $\ggo$ and identified as usual with the left invariant vector fields of $G$. The vector space generated by the vectors $v_1, \hdots v_j$ will be denoted $\la v_1, \hdots, v_j \ra$. Throughout this paper we assume that Lie groups are simply connected (and so there is a bijection between the set of simply connected Lie groups and the set of Lie algebras).

\section{Preliminaries}

Since we are interested on left invariant structures, our work is reduced to the Lie algebras of the corresponding Lie groups. Thus the  first step is to determine the different classes of Lie algebras, which are our topic of study. Let $\ggo$ be a four dimensional solvable real Lie algebra. The following table  shows the different types of such Lie algebras (see \cite{D}, \cite{Mu} or for a complete proof \cite{A-B-D-O}). 

\begin{prop} \label{clasi} Let $G$ be a simply connected solvable four dimensional real Lie group. Then $G$ has a Lie algebra $\ggo$ equivalent to one and only one of the Lie algebras listed below:

\

{\small
\begin{tabular}{|c|c|c|}\hline 
{\,\,\,{$\ggo$}\,\,\,} & {\qquad \rm{Lie bracket relations}\qquad \quad}&{\qquad $\ggo'$\qquad }  \\ \hline 
&  & \\
{\qquad $\aa_4$\qquad } & {$\qquad \qquad [e_i,e_j] = 0, \quad \mbox{\rm for all } e_i, e_j$ \qquad \qquad }& {\qquad $0$ \qquad}  \\ 
& & \\\hline 
& & \\
{$\rr\hh_3$} & {$[e_1,e_2] = e_3 $}  & { $\la e_3 \ra$ }\\
 {} & & \\ \hline 
\end{tabular}} 
\pagebreak

{\small
\begin{tabular}{|c|c|c|} \hline 
&& \\
{\qquad $\rr\rr_3$\qquad } & {\qquad $[e_1,e_2] =  e_2,\,[e_1,e_3] = e_2+ e_3 $\qquad }  & {\qquad$\la e_3, e_2 \ra$ \qquad}\\
  & & \\ \hline 
& & \\
{$\rr\rr_{3,\lambda}$} & {$[e_1,e_2] =e_2, [e_1,e_3]=\lambda e_3$} &{$\la e_2, \lambda e_3\ra$}\\
{$\lambda \in [-1,1]$} & & \\ \hline 
& & \\
{$\rr \rr'_{3,\gamma}$} & {$ [e_1, e_2] = \gamma e_2 -  e_3,[e_1, e_3] =  e_2 + \gamma e_3$} &{$ \la e_2,e_3\ra $}\\ 
{$\gamma \ge 0$} & {} &\\ \hline 
& & \\
{$\rr_2 \rr_2$} & {\qquad $[e_1, e_2] = e_2,\,[e_3,e_4] = e_4,\, \qquad  $} &{$\la e_2, e_4\ra$} \\ 
&  & \\\hline 
&  & \\
{$\rr_2'$} & {$[e_1,e_3] = e_3,\,[e_1,e_4] = e_4,$} &{$\la e_3,e_4\ra$}\\ 
&{$\, [e_2,e_3]= e_4, \,[e_2,e_4]=-e_3 $} & \\ \hline 
& & \\
{$\nn_4$} & {$[e_4,e_1] = e_2,\,[e_4, e_2] = e_3 $} &{$\la e_2, e_3 \ra$}\\
& &  \\ \hline &  & \\
{$\rr_4$} & {$[e_4,e_1]=e_1,\, [e_4,e_2] = e_1 + e_2, [e_4, e_3]=e_2 + e_3$}& {$\la e_1,e_2, e_3 \ra $} \\ 
&  &\\\hline 
&  & \\  
{$\rr_{4,\mu}$} & {$[e_4,e_1] = e_1,\, [e_4,e_2] = \mu e_2,\,[e_4,e_3] = e_2 + \mu e_3 $} & {$ \la e_1, e_2, \mu e_3 \ra$} \\  
{$\mu \in \RR$} & {} &  \\ \hline
&& \\
{$\rr_{4,\alpha, \beta}$}& & \\
{$\alpha \beta \ne 0$} & {$[e_4,e_1] = e_1,\, [e_4,e_2] = \alpha e_2,\, [e_4,e_3] = \beta e_3$} & {$ \la e_1, e_2,  e_3\ra $} \\
{$-1 \le \alpha \le \beta \le 1$} & & \\ \hline 
&  & \\ 
{$\rr'_{4,\gamma, \delta}$} & {$[e_4,e_1] = e_1,\, [e_4,e_2] = \gamma e_2 - \delta e_3, $} & {$\la e_1, e_2, e_3\ra $}  \\ 
{$\gamma \in \RR, \delta > 0$} & {$[e_4,e_3] = \delta e_2 +\gamma e_3$} & \\  \hline
& & \\ 
{$\dd_4$} & {\quad $[e_1,e_2]=e_3,\, [e_4,e_1] = e_1,\, [e_4,e_2] = -e_2$\quad } &{ $\la e_1, e_2, e_3\ra$}  \\
& & \\ \hline 
& & \\ 
{$\dd_{4,\lambda}$} & {$[e_1,e_2]=e_3,\, [e_4,e_3] = e_3,\,$} & {$\la e_1, e_2, e_3 \ra$} \\ 
{$\lambda \ge \frac12$} & {$ [e_4, e_1]=\lambda  e_1,\,\,[e_4, e_2]=(1-\lambda) e_2$}  & \\  \hline 
& & \\
{$\dd'_{4,\delta}$} & {$[e_1,e_2]=e_3,\, [e_4, e_1]=\frac{\delta}2  e_1- e_2,$}& {$\la e_1, e_2, e_3 \ra$}\\ 
{$\delta \ge 0$} & {$[e_4,e_3] = \delta e_3,\,\,[e_4, e_2]=e_1+ \frac{\delta}2 e_2$} & \\  \hline
& & \\
{$\hh_4$} & { $[e_1,e_2]=e_3,\, [e_4,e_3] = e_3,\,\,$}  & {$\la e_1, e_2, e_3 \ra,$}\\ 
{} &{$ \qquad[e_4, e_1]=\frac{1}2 e_1,\,\,[e_4, e_2]= e_1 + \frac{1}2 e_2 $\qquad } & \\  \hline 
\end{tabular}} 

\begin{center} \label{clases}
\rm{Table} \ref{clases} 
\end{center} 
\end{prop}

\begin{remark} For an explanation concerning the Lie groups which admit a cocompact quotient, see for example the work of Oprea and Tralle \cite{O-T}.
\end{remark}

\begin{remark}
The case $\rr_2\rr_2$ corresponds to the Lie algebra $\aff(\RR) \times \aff(\RR)$, the case $\rr_2'$ corresponds to the Lie algebra $\aff(\CC)$, $\rr\rr_{3,-1}$ is the trivial extension of $\ee(1,1)$, the Lie algebra  corresponding to the Lie group of rigid motions of the Minkowski 2-space;
$\rr'_{3,0}$ is the trivial extension of $\ee(2)$ the Lie algebra of the  Lie group of rigid motions of $\RR^2$; $\rr\hh_3$ is the trivial extension of
the three-dimensional Heisenberg Lie algebra  denoted by $\hh_3$. 
\end{remark}

\section{Complex structures}

An {\it invariant complex structure} on a real Lie group $G$ is a complex structure on the underlying 
manifold such that left multiplication $\ell_g,\, g \in G$ (but not necessarily right multiplication) by elements of the group are holomorphic with respect to this structure. Because of the Newlander-Niremberg theorem \cite{N-N} we have the following equivalent definition:

\begin{defn}\label{def1} An invariant complex structure on a real Lie group $G$ is an endomorphism $J$ of $\ggo$ the Lie algebra of $G$, such that: 
 
(i) $J^2 = - Id$; 
 
(ii) $0 = [JX,JY] - [X,Y] - J[JX,Y] - J[X,JY] \qquad \forall X, Y \in \ggo$ 
 
\end{defn} 
 
Condition ii) is called the {\it integrability condition} of $J$ and usually for $X, X$ in $\ggo$ one denotes by $N_J(X,Y)$ the right hand of equality ii) and thus the integrability condition is announced as $N_J(X,Y)=0$ for all $X,Y \in \ggo$. 
 
By extending $J$ in a natural way to $\ggo^{\CC} =\ggo \otimes_{\RR} \CC$ the complexification of the real Lie algebra $\ggo$, it is possible to have an equivalent condition to that given in Definition \ref{def1}. We denote by $\sigma$ the conjugation in $\ggo^{\CC}$ with respect to the real form $\ggo$, that is, $\sigma (X + i Y) = X - i Y$, $X, Y \in \ggo$. 
 
\begin{prop}\label{prop2} A real Lie group $G$ has a left invariant complex structure if and only if $\ggo^{\CC}$ admits a decomposition as a direct sum of vector spaces: 
\begin{equation} 
\ggo^{\CC} = \qq \oplus \sigma \qq \label{def2} 
\end{equation} 
where $\qq$ is a complex subalgebra of $\ggo^{\CC}$. 
\end{prop} 
\begin{proof} Condition (i) in Definition \ref{def1} gives a decomposition of $ \ggo^{\CC}$ into a direct sum of subspaces (the eigenspaces of $J$).  Condition (ii) is equivalent to the fact that these subspaces are subalgebras. Conversely, given a complex subalgebra $\qq$ satisfying (\ref{def2}), let $J$ be the almost complex structure defined on $\ggo^{\CC}$ by $JX = -i X$, $J\sigma X = i \sigma X$, for $X \in \qq$. Since $J \sigma = \sigma J$, it is possible to define $J$ on $\ggo$. Using the fact that $\qq$ is a subalgebra, it is not hard to see that $J$ is integrable.  
 
\vskip .3cm 
 
Thus, the previous proposition shows that there exists a one to one correspondence between left invariant complex structures $J$ and subalgebras $\qq$ satisfying  (\ref{def2}). 
\end{proof} 
 
The subalgebras which satisfy (\ref{def2}) are called {\it (invariant) complex subalgebras}. 

We say that two invariant complex structures  $J_1$ and $J_2$ on a real Lie group are {\it equivalent} if there exists $x \in \Aut(\ggo)$ such that $x J_1 = J_2 x$. In terms of complex subalgebras, we say that two invariant complex subalgebras, $\qq_1$ and $\qq_2$ are equivalent if there exists $y \in \Aut(\ggo^{\CC})$ such that $y \sigma = \sigma y$ and $y \qq_1 = \qq_2$. Note that both equivalent conditions are equivalent.

As we said in the Introduction invariant complex structures were classified by J. Snow \cite{Sn} and G. Ovando \cite{O}. They search for non equivalent classes and determine the underlying complex manifold related to the complex structure. But they did not investigate the space of complex structures. It is our aim now to give the general form of any complex subalgebra for a simply connected four-dimensional solvable real Lie algebra $\ggo$. Observe that only partial results can be found in the papers of \cite{Sn} \cite{O}, so the third  column ($\mathcal Q$) of the table of the following Proposition completes this information.

\begin{prop} \label{complejas} Let $G$ be a simply connected  four dimensional solvable real Lie group. Then in the notation of Proposition \ref{clasi}, the following table shows the non-equivalent classes of complex subalgebras $\qq$ (\cite{Sn}, \cite{O}),and the general form of any such complex subalgebra, denoted by $\mathcal Q$, only in those cases when they exist:

\begin{center} 
\small{
\begin{tabular}{|c|c|c|}\hline 
{{\rm Case}} & {{$\qq$}} & {${\mathcal Q}$}\\ \hline 
&& \\
{$\rr\hh_3$} & {$\langle e_1+ie_2,e_3+ie_4\rangle $} & {$\la   e_1 + b_1e_2 + d_1 e_4, e_3 + d_2e_4\ra$} \\ 
{}&& {$\Im b_1 \Im d_2 \ne0$} \\
\hline 
&& \\
{$\rr\rr_{3,0}$} & {$\langle e_1+ie_2,e_3+ie_4\rangle $} & {$\la   e_1 + b_1e_2 + d_1 e_4, e_3 + d_2e_4\ra$} \\ 
{}&& {$\Im b_1 \Im d_2 \ne0$} \\
\hline 
&& \\
{$\rr\rr_{3,1}$} & {$\langle e_1+ie_4,e_3+ie_2\rangle $} & {$\la   e_1 + b_1e_2 + d_1 e_4, b_2e_2 + e_3\ra$} \\ 
{}&& {$\Im b_2 \Im d_1 \ne0$} \\
\hline 
&& \\
{$\rr\rr'_{3,0}$} & {$\langle e_1 + ie_4, e_2+i e_3\ra$} & {} \\ 
&& {$ \la e_1 + c_1e_3 + d_1 e_4, e_2 - (2 \gamma +\varepsilon i)e_3\ra$}\\ \cline{1-2} 
&& \\
{$\rr\rr'_{3,\gamma}$} & {$\langle e_1 + ie_4, e_2-( 2\gamma + i) e_3\ra, $} & { $\varepsilon = \pm 1,\, \Im d_1 \ne 0$} \\
{$\gamma \ne 0$}&{$\la e_1+ie_4, e_2-(2\gamma -i)e_3 \ra $} & {} \\
&&\\\hline 
&& \\
{$\rr_2\rr_2$} & {$\langle e_1 +ie_2, e_3+ie_4\rangle$} & {$\la e_1 + b_1e_2+d_1 e_4, e_3+d_2e_4\ra$} \\ && {$ \Im b_1 \Im d_2 \ne 0$}\\
\hline 
&& \\
 &  &  {$\la e_1 + c_1 e_3 + d_1 e_4, e_2 -d_1 e_3+c_1 e_4\ra$}\\
 & {$\langle e_1+ie_3, e_2+ie_4\rangle$}& {$(\Im c_1)^2 +(\Im d_1)^2\ne 0$} \\ & & \\
{$\rr_2'$}&  &  {$\la e_1 + \varepsilon i e_2 + d_1 e_4,  e_3- \varepsilon i e_4\ra,\, \varepsilon = \pm 1$}\\
 &{$\langle e_1+ib_1e_2, e_3+ie_4\rangle,$} & {} \\  
& & {$\la e_1 + b_1e_2 + d_1 e_4,  e_3 + \varepsilon i e_4\ra\,$}\\
&{$\, \Im b_1 \ne 0 $}& {$\varepsilon =\pm 1,\, \Im b_1 \ne 0 $}\\
 &{}  & {}\\ \hline 
& & \\
{$\rr_{4,1}$} & {$\langle e_4+ie_3, e_1+ie_2\rangle$} & {$\la e_4+ b_1 e_2+c_1 e_3, e_1 + b_2 e_2\ra$} \\ & &{$ \Im c_1 \Im b_2 \ne 0$} \\ \hline
&&\\
{$\rr_{4,\alpha, 1}$} & {$\langle e_4+ie_2,e_1+ie_3\rangle$} & {$\la e_4+ b_1 e_2+c_1 e_3, e_1 + c_2 e_3\ra$} \\{$\alpha \ne 1$} &&{$\Im b_1 \Im c_2 \ne 0$} \\\hline 
&& \\
{$\rr_{4,\alpha, \alpha}$} & {$\langle e_4+ie_1,e_2+ie_3\rangle$} & {$\la e_4 + a_1e_1+ c_1 e_3, a_2e_1 + e_2 + c_2 e_3\ra$} \\ 
{$\alpha \ne 0,1$} & {} & {$\Im c_1 \Im a_2-\Im a_1 \Im c_2 \ne 0$} \\ \hline  
& & \\
{$\rr'_{4,\gamma, \delta}$} & {$\la e_4+ie_1, e_2+ie_3 \ra, $} & {$ \la e_4+ a_1 e_1+c_1 e_3, e_2 + \varepsilon i  e_3\ra$} \\
{$\,$} & {$\la e_4+ie_1, e_2-ie_3\ra$} & {$\varepsilon = \pm 1,\, \Im a_1 \ne 0$} \\ && \\ \hline 
& & {}\\
{$\dd_4$} & {$\la e_1-ie_3,e_4+ie_2\ra, $} & {$\la e_4+ b_1 e_2 + c_1 e_3, e_1 - b_1 e_3\ra \, \Im b_1 \ne 0$} \\ 
 &{$\la e_4 +i( e_2+e_3), e_1 - i e_3\ra$} & {$\la e_4+ a_1 e_1+c_1 e_3, e_2 -a_1 e_3\ra\, \Im a_1 \ne 0$}\\ && \\\hline
\end{tabular}} 
 \end{center} 

\begin{center}
\small{
\begin{tabular}{|c|c|c|}\hline
&& \\
{$\dd_{4,1}$} & { $\langle e_4 - ie_1, e_2 + i e_3\rangle$} &  {$\la e_4 - c_2e_1 + c_1 e_3, e_2 + c_2 e_3\ra $}  \\ && {$ \Im c_2 \ne 0$}\\
\hline  
& & \\
{} & {$\langle e_4+ie_3, e_1+ie_2 \rangle$} &  {$\la e_4+ b_1 e_2+ c_1 e_3,  e_1 + b_2 e_2 + 2 b_1 e_3  \ra$} \\
 & & {$ 2\Im ^2  b_1 - \Im c_1 \Im b_2 \ne 0$}\\  
{$\dd_{4,1/2}$}& {$\langle e_4+ie_3, e_1 - ie_2\rangle$} & \\
 &  & {$\la e_4+ a_1 e_1+ c_1 e_3, a_2 e_1 + e_2-2 a_1 e_3 \ra$} \\
 & {$\langle e_4+ie_1, e_2-2ie_3\rangle$}& {$ 2\Im ^2a_1+\Im c_1 \Im a_2  \ne 0$}\\ 
 {} &  & \\ \hline
 && {}\\
{$\dd_{4,\lambda}$} & {$ \langle e_4- \lambda ie_1, e_2+i e_3\rangle$} & {$\la e_4-\lambda c_2 e_1+ c_1 e_3, e_2 + c_2 e_3 \ra\, $}  \\ 
{$ \lambda \ne 1/2,1$}& {$\langle e_4+(1 - \lambda) ie_2, e_1+ i e_3 \rangle$} & {$\la e_4+ (1 - \lambda)c_2 e_2+c_1 e_3, e_1+c_2 e_3 \ra\,$}\\ && {$\Im c_2 \ne 0 $} \\\hline
& {} & \\
{$\dd'_{4,0}$} & {$\langle e_4-ie_3,e_1+ie_2\rangle$} & {$\la e_4+ c_1 e_3-\varepsilon i c_2 e_2, c_2 e_3 + e_1 + \varepsilon i e_2\ra$} \\ 
 & {$\langle e_4+ie_3,e_1+ie_2 \rangle$} & {$\la e_4+ c_1 e_3-\varepsilon ic_2 e_2, c_2 e_3 + e_1 + \varepsilon i e_2\ra$} \\ &&{$ \varepsilon = \pm 1,\, \Re c_2 \Im c_2 + \Im c_1 \ne 0$}\\
\hline
& & {}\\ 
{$\dd'_{4,\delta}$} & {$ \la e_4+ie_3,e_1-ie_2 \ra$} & {}{$\la e_4+ c_1 e_3+ (\delta/2 -\varepsilon i)c_2 e_2 ,  $} \\ 
{$\delta > 0$} & {$\la e_4-ie_3, e_1+ie_2 \ra $} &{$ e_1 +\varepsilon i e_2+ c_2 e_3\ra,$} {$\,\, \varepsilon =\pm1,\,$}\\ 
 &{$\la e_4+ie_3, e_1+ie_2 \ra,\,$} &  \\
 & {$\la e_4-ie_3, e_1-ie_3 \ra$} & {$ \varepsilon \Im c_1 + \Im c_2(-\frac{\delta}2\Im c_2 - \varepsilon \Re c_2) \ne 0$}\\ && \\ \hline 
& & \\
{$\hh_4$} & {$\langle e_4+ie_2,e_1+2ie_3\rangle$} & {$\la e_4+ b_1 e_2+  c_1 e_3,  2b_1 e_3 + e_1 \ra$}\\ 
{} & {} & {$\Im b_1 \ne 0$} \\  \hline
\end{tabular} }
\end{center} 
\begin{center} 
{\rm Table} \ref{complejas} 
\end{center} 
\end{prop}

\begin{proof} It is not our purpose to give a complete proof, which consists of many technical computations, but to write the idea of the proof. In fact the mentioned papers \cite{Sn} \cite{O} (or in the three dimensional commutator case also  \cite{O1}) show explicitly the computations in particular cases.  These computations are repeated in each case and follow from the following reasoning.

If there exists  $\qq$ an invariant complex subalgebra of $\ggo^{\CC}$, then there are elements $U, V \in \ggo^{\CC}$ of the form:
$$
U = e_4 + a_1 e_1 + b_1 e_2 + c_1 e_3, \qquad V = a_2 e_1 + b_2 e_2 + c_2 e_3
$$
(by changing subindices if necessary) which are a basis of $\qq$ and
\begin{equation}
[ U, V ] = \beta V \qquad \qquad \beta \in \CC\label{ecu1}
\end{equation}
 
From (\ref{ecu1}) we get a system of equations, such that, the existence of solutions (that is, the existence of coefficients $a_i, b_i, c_i$ and $\beta$) is equivalent to the existence of complex subalgebras,when we imposed an extra condition to make of the set  $\{ U, V, \sigma U , \sigma V\}$ a basis of $\ggo^{\ast}$.
\end{proof}

\begin{remark} For explicit computations of the authomorphism group of these Lie algebras see \cite{O2} or in the three dimensional commutator case see \cite{O1}.
\end{remark}

\begin{remark} 
A complex structure $J$ is said to be {\it abelian} if it satisfies $[JX, JY] = [X,Y]$ for all $X, Y \in \ggo$, or equivalent in terms of complex subalgebras, if $\qq$ (or $\sigma \qq$) is abelian.

In dimension four only the solvable real Lie algebras with a commutator of dimension less than three can admit a abelian complex structure. In fact, that occurs in the cases $\aa_4$, $\rr_{4,0,0}$,  $\rr\hh_{3}$, $\rr_2\rr_2$, $\rr'_2$, $\dd_{4,1}$. 

\vspace{.2cm}

A complex structure $J$  introduces on $\ggo$ a structure of complex algebra if $J \circ ad_X = ad \circ J X$ for all $X \in \ggo$, and so $(\ggo, J)$ is a {\it complex algebra}, and that means that the corresponding simply connected Lie group is also complex, that is, left and right multiplication by elements of the Lie group are holomorphic maps. Only the Lie algebras $\aa_4$, the abelian one, and $\rr'_2$, $\aff(\CC)$, admit a such kind of complex structures. This fact can be proved by checking the Table \ref{complejas}.

\vspace{.2cm}

An hypercomplex structure on a Lie group $G$ is a pair of complex structures $J_1, J_2$ such that $J_1J_2 = -J_2 J_1$. The hypercomplex structures on four dimensional real Lie algebras were classified by L. Barberis (\cite{B}).
\end{remark}

\section{symplectic structures}

Let $M$ be a differential manifold, denote as usual by $\Lambda^{k}(M)$ the space of all smooth $k$-forms, and by $\Lambda^{\ast}(M)$ the algebra of all differential forms. Denote by $d: \Lambda^{\ast}(M) \to \Lambda^{\ast}(M)$ the unique antiderivation of degree +1 which satisfies

i) $d^2= 0$

ii) Whenever $ f \in C^{\infty}(M) = \Lambda^o(M)$, $df$ is the differential of $f$.

\

Let us recall some properties of $d$:

iii) $d(\omega_1 \wedge \omega_2) = d\omega _1 \wedge \omega_2 + (-1)^r \omega _1 \wedge d\omega_2 \quad \omega_1 \in \Lambda^r(M),\, \omega_2 \in \Lambda^{\ast}(M)$

iv) $$\begin{array}{rcl}
d\omega(Y_0, \hdots Y_p) & = & \sum_{i=0}^p (-1)^i Y_i \omega(Y_0, \hdots, Y_i, \hdots , Y_p) + \\
& & + \sum_{i<j}(-1)^{i+j}\omega([Y_i,Y_j], \hdots, Y_i, \hdots, Y_j , \hdots Y_p)\end{array}$$

In particular, when $M=G$ is a Lie group and the considered forms are left invariant, then the first sum of the right hand in iv) vanishes.

\begin{defn} A symplectic structure on a 2n-dimensional manifold $M$ is a  2-form $\omega \in \Lambda^2(M)$ such that:

i) $\omega$ is closed , that is $d \omega =0$,

ii) $\omega$ is a volume form on $M$.

\end{defn}

Our aim now is to study the existence of left invariant symplectic structures on a four dimensional solvable real Lie group. Left invariant symplectic structures are those which are invariant by left translation of elements of the group, that is,
$$
 \ell_g^{\ast} \omega = \omega \circ {\ell_g}_ {\ast}
$$
The left invariance property allows to work on the Lie algebra $\ggo$. By denoting by $\{e^i\}$ the dual basis on $\ggo^{\ast}$ of the basis $\{e_i\}$ on $\ggo$, the following proposition  summerizes the results in the notation compatible with Proposition (\ref{clasi}).

\begin{prop} \label{simplecticas} Let $G$ be a solvable real Lie group of dimension four. Then the following table listes such cases where there exists a left  invariant symplectic structure $\omega$ on $G$, give the general form of such one by imposing some extra condition:

\begin{center} 
\small{
\begin{tabular}{|c|c|c|}\hline 
{{\rm Case}} & {{$\omega$}} & {\rm Condition}\\ \hline 
&& \\
{$\rr\hh_{3}$} & {$a_{12}e^1 \wedge e^2+a_{13}e^1 \wedge e^3 + a_{14}e^1 \wedge e^4 + $} & {$a_{14}a_{23}-a_{13}a_{24}\ne 0$} \\ 
{}&{$a_{23}e^2 \wedge e^3+a_{24}e^2 \wedge e^4$}& \\
\hline &&\\
{$\rr \rr_{3,0}$}& {$a_{12}e^1 \wedge e^2 + a_{13}e^1 \wedge e^3 +a_{14} e^1 \wedge e^4 +a_{34} e^3 \wedge e^4$}& {$a_{12} a_{34} \ne 0$} \\
&&\\ \hline
&&\\
{$\rr \rr_{3,-1}$}& {$a_{12}e^1 \wedge e^2 + a_{13}e^1 \wedge e^3 +a_{14} e^1 \wedge e^4 +a_{23} e^2 \wedge e^3$}& {$a_{14} a_{23} \ne 0$} \\
&&\\ \hline
&&\\
{$\rr\rr'_{3,0}$} & {$a_{12}e^1 \wedge e^2 + a_{13}e^1 \wedge e^3+a_{14}e^1 \wedge e^4 + a_{23}e^2 \wedge e^3$} & {$a_{14}a_{23}\ne 0$} \\ && \\ \hline
&& \\
{$\rr_2\rr_2$} & {$a_{12}e^1 \wedge e^2 + a_{13}e^1 \wedge e^3 + a_{34}e^3 \wedge e^4$} & {$a_{12}a_{34}\ne 0$} \\ && \\
\hline 
&& \\
{$\rr'_2$} & {$a_{12}e^1 \wedge e^2 + a_{13-24}(e^1 \wedge e^3 - e^2 \wedge e^4) +  $} & {$a_{14+23}^2 + a_{13-24}^2\ne 0$} \\ & {$+a_{14+23}(e^1 \wedge e^4+e^2 \wedge e^3)$} & \\ \hline 
&& \\
{$\nn_{4}$} & {$a_{12}e^1 \wedge e^2 + a_{14}e^1 \wedge e^4+a_{24}e^2 \wedge e^4 + a_{34}e^3 \wedge e^4$} & {$a_{12}a_{34}\ne 0$} \\ && \\ \hline
& & \\
{$\rr_{4,0}$} & {$a_{14}e^1 \wedge e^4 + a_{23}e^2 \wedge e^3+a_{24}e^2 \wedge e^4 + a_{34}e^3 \wedge e^4$} & {$a_{14}a_{23}\ne 0$} \\ && \\ \hline
& & \\
{$\rr_{4,-1}$} & {$a_{13}e^1 \wedge e^3 + a_{14}e^1 \wedge e^4+a_{24}e^2 \wedge e^4 +a_{34}e^3 \wedge e^4$} & {$a_{13}a_{24}\ne 0$}\\ 
{} & {} & {} \\  \hline
&&\\
{$\rr_{4,-1, \beta}$} & {$a_{12}e^1 \wedge e^2 + a_{14}e^1 \wedge e^3+a_{24}e^2 \wedge e^4 ++ a_{34}e^3 \wedge e^4$}  & {$a_{14}a_{23}\ne 0$} \\ {$\beta \ne -1,0,1$} & {}& {}
\\\hline{} && \\
{$\rr_{4,-1, -1}$} & {$a_{12}e^1 \wedge e^2 + a_{13}e^1 \wedge e^3+a_{14}e^1 \wedge e^4+ $} & {$a_{12}a_{34} - a_{13}a_{24}\ne 0$} \\
{} &{$+ a_{24}e^2 \wedge e^4+a_{34}e^3 \wedge e^4$}& \\ \hline
&&\\
{$\rr_{4,-1, 1}$} & {$a_{12}e^1 \wedge e^2 + a_{14}e^1 \wedge e^4+a_{23}e^2 \wedge e^3 + $} & {$a_{12}a_{34}+a_{14}a_{23}\ne 0$} \\{} &{$+a_{24}e^2 \wedge e^4+a_{34} e^3 \wedge e^4$}&  \\ \hline
&&\\
{$\rr_{4,\alpha, -\alpha}$} & {$a_{14}e^1 \wedge e^4 + a_{23}e^2 \wedge e^3+a_{24}e^2 \wedge e^4 +a_{34}e^3 \wedge e^4$} & {$a_{14}a_{23}\ne 0$} \\{$\alpha \ne -1,0$} &{}& \\ \hline
&&\\
{$\rr'_{4,0, \delta}$} & {$a_{14}e^1 \wedge e^4 + a_{23}e^2 \wedge e^3+a_{24}e^1 \wedge e^4 +$} & {$a_{14}a_{23}\ne 0$}\\ 
{$\delta \ne 0$} & {$+ a_{34} e^3 \wedge e^4$} & {} \\  \hline
  
&& \\
{$\dd_{4,1}$} & { $a_{12-34}(e^1 \wedge e^2 - e^3 \wedge e^4) +a_{14}e^1 \wedge e^4 + a_{24}e^2 \wedge e^4$} & {$a_{12-34}\ne 0$} \\ && \\ \hline
&  & \\
{$\dd_{4,2}$} & {$a_{12-34} (e^1\wedge e^2-e^3 \wedge e^4)+a_{14}e^1 \wedge e^4 +$} & {$a_{12-34}^2+a_{14}a_{23}\ne 0$} \\
 &{$+ a_{23} e^2 \wedge e^3+ a_{24} e^2 \wedge e^4$} & \\ \hline
&&\\
{$\dd_{4,\lambda}$} & {$a_{12-34} (e^1\wedge e^2-e^3 \wedge e^4)+a_{14}e^1 \wedge e^4 +$} & {$a_{12-34}\ne 0$} \\
 {$\lambda  \ne 1,2$}& {$+ a_{24} e^2 \wedge e^4$}& \\ \hline
\end{tabular}} 
 \end{center}
 
\pagebreak

\begin{center}
\small{
\begin{tabular}{|c|c|c|}\hline 
&&\\
{$\dd'_{4,\delta}$} & {$a_{-12+\delta34} (-e^1\wedge e^2+\delta e^3 \wedge e^4)+a_{14}e^1 \wedge e^4 +$} & {$a_{-12+\delta 34}\ne 0$} \\
 {$\delta  \ne 0$}& {$+ a_{24} e^2 \wedge e^4$}& \\ \hline
&&\\
{$\hh_4$} & {$a_{12-34} (e^1\wedge e^2-e^3 \wedge e^4)+a_{14}e^1 \wedge e^4 +$} & {$a_{12-34}\ne 0$} \\
&{$+ a_{24} e^2 \wedge e^4$}&\\\hline
\end{tabular} }
\end{center}
\begin{center} 
{\rm Table} \ref{simplecticas} 
\end{center} 
\end{prop}

\begin{proof} As all cases should be worked out in a similar way, we will give as example the proof of the case $\rr'_2$, that is the Lie algebra which corresponds to $\aff (\CC)$. 

From the definition of the Lie bracket it follows that $d e^1 =0=de^2$ and simple computations  following the properties mentioned at the beginning of the section allow to get $-de^3 = e^1\wedge e^3-e^2\wedge e^4$, and $-de^4 = e^1 \wedge e^4 + e^2 \wedge e^3$. At the next level we have:
$$ d(e^1 \wedge e^3) = -e^1\wedge e^2 \wedge e^4,\, d(e^1 \wedge e^4) = e^1\wedge e^2 \wedge e^3,\, d(e^2 \wedge e^3) = -e^1\wedge e^2 \wedge e^3$$
$$d(e^2 \wedge e^4) =- e^1\wedge e^2 \wedge e^4,\qquad d(e^3 \wedge e^4) = -2e^1\wedge e^3 \wedge e^4$$
Thus any 2-form $\omega$ that is closed has the form $\omega = a_{12}e^1 \wedge e^2+a_{13-24}(e^1\wedge e^3- e^2 \wedge e^4) +a_{14+23} (e^1\wedge e^4+e^2 \wedge e^3)$. Now $\omega$ is symplectic if it satisfies the conditions mentioned at the Table (\ref{simplecticas}) for $\rr'_2$.
\end{proof}

\begin{cor} A four dimensional solvable Lie algebra admits an exact symplectic structure if and only if $\ggo$ is one of the following attached with the respective symplectic structure

\begin{center}
\small{
\begin{tabular}{|c|c|c|}\hline 
{{\rm Case}} & {{$\omega$}} & {\rm Condition}\\ \hline 
{$\rr_2\rr_2$} & {$a_{12}e^1 \wedge e^2 + a_{34}e^3 \wedge e^4$} & {$a_{12}a_{34}\ne 0$} \\\hline 
{$\rr'_2$} & {$ a_{13-24}(e^1 \wedge e^3 - e^2 \wedge e^4) +a_{14+23}(e^1 \wedge e^4+e^2 \wedge e^3)  $} & {$a_{14+23}^2 + a_{13-24}^2\ne 0$} \\ \hline 
{$\dd_{4,1}$} & { $a_{12-34}(e^1 \wedge e^2 - e^3 \wedge e^4) +a_{14}e^1 \wedge e^4$} & {$a_{12-34}\ne 0$} \\ \hline
{$\dd_{4,\lambda}\lambda \ne 1$} & { $a_{12-34}(e^1 \wedge e^2 - e^3 \wedge e^4) +a_{14}e^1 \wedge e^4 + a_{24}e^2 \wedge e^4$} & {$a_{12-34}\ne 0$} \\ \hline
{$\dd'_{4,\delta}\delta \ne 0$} & {$a_{-12+\delta34} (-e^1\wedge e^2+\delta e^3 \wedge e^4)+a_{14}e^1 \wedge e^4 +a_{24}e^2 \wedge e^4$} & {$a_{-12+\delta 34}\ne 0$} \\\hline
{$\hh_4$} & {$a_{12-34} (e^1\wedge e^2-e^3 \wedge e^4)+a_{14}e^1 \wedge e^4 + a_{24} e^2 \wedge e^4$} & {$a_{12-34}\ne 0$} \\\hline
\end{tabular} }
\end{center}
\end{cor}

\begin{remark} Compare with the results obtained in \cite{Ca}. \end{remark}

\

\subsection{Appendix: Cohomology} To conclude this section we would like to compute the cohomology of these Lie algebras.
 Recall that for $i \ge 0$ the cohomology of the Lie algebra $\ggo$ over $\RR$ is given by the groups $H^i(\ggo)$, where
$$ H^i(\ggo) = \frac{\ker d: \Lambda^{i}(\ggo) \to \Lambda^{i+1}(\ggo)}{\Rg d:\Lambda^{i-1}(\ggo) \to \Lambda^{i}(\ggo)}$$
and $\ker d$ denotes the kernel of $d$ and $\Rg d$ denotes the image of $d$. The class of any element $\theta \in \Lambda (\ggo)$ will be denoted by $[\theta]$.

\begin{prop} \label{coho} The cohomology over $\RR$ of any four dimensional solvable real Lie algebra is presented in the following Table:

\begin{center} 
\small{
\begin{tabular}{|c|c|c|c|}\hline 
{{\rm Case}} & {{$H^1(\ggo)$}} & {${H^2(\ggo)}$} & {$H^3(\ggo)$}\\ \hline 
{$\rr \hh_3$} & {$[e^1][e^2][e^4]$} & {$[e^1 \wedge e^3][e^1 \wedge e^4]$} & {$[e^1 \wedge e^2 \wedge e^3] [e^1 \wedge e^3 \wedge e^4] $} \\ 
 & & {$[e^2 \wedge e^3][e^2 \wedge e^4]$} & {$[e^2 \wedge e^3 \wedge e^4]$} \\ \hline
{$\rr\rr_{3}$} & {$[e^1][e^4]$} & {$[e^1 \wedge e^4]$} & {$0$} \\
\hline
{$\rr \rr_{3,0}$} & {$[e^1][e^4]$} & {$[e^1 \wedge e^4]$} & {$[e^1 \wedge e^2 \wedge e^3]$} \\
\hline
{$\rr \rr_{3,-1}$} & {$[e^1][e^4]$} & {$[e^1 \wedge e^4][e^2 \wedge e^3]$} & {$[e^1 \wedge e^2 \wedge e^3][e^2 \wedge e^3 \wedge e^4]$} \\
\hline
{$\rr \rr_{3,\lambda},\,\lambda \ne 0,-1$} & {$[e^1][e^4]$} & {$[e^1 \wedge e^4]$} & {$0$} \\ 
\hline
{$\rr \rr'_{3,0}$} & {$[e^1][e^4]$} & {$[e^1 \wedge e^4][e^2 \wedge e^3]$} & {$[e^1 \wedge e^2 \wedge e^3][e^2\wedge e^3\wedge e^4] $} \\ 
\hline
{$\rr \rr'_{3,\gamma}$, $\gamma \ne 0$} & {$[e^1][e^4]$} & {$[e^1 \wedge e^4]$} & {$0$} \\ 
\hline
{$\rr_2 \rr_2$} & {$[e^1][e^3]$} & {$[e^1 \wedge e^3]$} & {$0$} \\ 
\hline
{$\rr_2'$} & {$[e^1][e^2]$} & {$[e^1 \wedge e^2]$} & {$0$} \\ 
\hline
{$\nn_4$} & {$[e^1][e^4]$} & {$[e^1 \wedge e^2][e^3 \wedge e^4]$} & {$[e^1 \wedge e^2 \wedge e^3] [e^2 \wedge e^3 \wedge e^4] $} \\
\hline
{$\rr_4$} & {$[e^4]$} & {$0$} & {$0$} \\
\hline
{$\rr_{4,0}$} & {$[e^3][e^4]$} & {$[e^2 \wedge e^3][e^2 \wedge e^4]$} & {$ [e^2 \wedge e^3 \wedge e^4]$} \\
\hline
{$\rr_{4,-1}$} & {$[e^4]$} & {$[e^2 \wedge e^4]$} & {$0$} \\ 
\hline
{$\rr_{4,-1/2}$} & {$[e^4]$} & {$0$} & {$[e^1 \wedge e^2 \wedge e^3]$} \\ 
\hline
{$\rr_{4,\mu},\, \mu\ne -1,-1/2,0$} & {$[e^4]$} & {$0$} & {$0$} \\ 
\hline
{$\rr_{4,-1, -1}$} & {$[e^4]$} & {$[e^1 \wedge e^3]$} & {$[e^1 \wedge e^2 \wedge e^4] [e^1\wedge e^3\wedge e^4] $} \\ 
\hline 
{$\rr_{4,-1, \beta}$} {$\beta \ne -1,-\alpha$}  
& {$[e^4]$} & {$[e^1 \wedge e^2]$} & {$[e^1 \wedge e^2 \wedge e^4]  $} \\
\hline 
 

{$\rr_{4,-1, 1}$} & {$[e^4]$} & {$[e^1 \wedge e^2][e^2 \wedge e^3]$} & {$[e^1 \wedge e^2 \wedge e^4] [e^2 \wedge e^3 \wedge e^4]$} \\
 \hline
{$\rr_{4,\alpha, -\alpha}$} {$\alpha \ne -1,0$}
& {$[e^4]$} & {$[e^2 \wedge e^3]$} & {$[e^2 \wedge e^3 \wedge e^4]  $} \\\hline
{$\rr_{4,\alpha, \beta}\,$}& & {} & \\
{$\alpha \ne -1,0,-\beta $},{$\alpha - \beta = -1$}
 & {$[e^4]$} & {$0$} & {$[e^1\wedge e^2\wedge e^3]  $} \\ \hline
{$\rr_{4,\alpha, \beta}, \,\,\alpha, \beta$}& & {} & \\
{\rm not as above }  {$\alpha - \beta \ne -1$}& {$[e^1]$} & {$0$} & {$0$} \\ \hline
{$\rr'_{4,0, \delta}$} {$\delta \ne 0$}& {$[e^4]$} & {$[e^2 \wedge e^3]$} & {$ [e^2 \wedge e^3 \wedge e^4]$} \\ 
\hline 
{$\rr'_{4,-1/2, \delta}$} {$\delta \ne 0$}& {$[e^4]$} & {$0$} & {$ [e^1 \wedge e^2 \wedge e^3] $} \\ 
\hline
{$\rr'_{4,\gamma, \delta}$} {$\gamma \ne -1/2, 0, \delta \ne 0$}& {$[e^4]$} & {$0$} & {$0$} \\ 
\hline
{$\dd_4$} & {$[e^4]$} & {$0$} & {$[e^1 \wedge e^2 \wedge e^3]$} \\
\hline
{$\dd_{4,1}$} & {$[e^2][e^4]$} & {$ [e^2 \wedge e^4]$} & {$0$} \\ 
\hline
{$\dd_{4,2}$} & {$[e^4]$} & {$[e^2 \wedge e^3]$} & {$[e^2 \wedge e^3 \wedge e^4] $} \\ 
\hline
{$\dd_{4,\lambda}$}  {$ \lambda  \ne 1,2\quad $}& {$ [e^4]$} & {$0$} & {$0$} \\
\hline
{$\dd'_{4,0}$} & {$[e^4]$} & {$0$} & {$[e^1 \wedge e^2 \wedge e^3][e^1\wedge e^2 \wedge e^4]$} \\
\hline 
{$\dd'_{4,\delta}$}{$\delta \ne 0$} & {$[e^4]$} & {$0$} & {$0$} \\ 
\hline
{$\hh_4$} & {$[e^4]$} & {$0$} & {$0$} \\
\hline
\end{tabular} }
\end{center}
\begin{center} 
{\rm Table} \ref{coho} 
\end{center} 
\end{prop}
\begin{proof} The cohomology can be obtained parallel to the computations made to get the Table (\ref{simplecticas}). Continuing with the case $\rr_2'$, (worked out in the previous Proposition (\ref{simplecticas})), in the proof of this Proposition we can see that $H^1(\ggo) =\{ [e^1], [e^2]\}$. From the computations at the next level it is possible to prove that $\theta \in H^2(\ggo)$ if and only if $\theta$ belongs to the class $[e^1 \wedge e^2]$. Since $e^1 \wedge e^2 \wedge e^3, e^1 \wedge e^2 \wedge e^4$ and $e^1 \wedge e^3 \wedge e^4$ are in the image of $d:\Lambda^{2}(\ggo) \to \Lambda^{3}(\ggo)$, to get $H^3(\ggo)$ we need to compute extra only the following: $d(e^2 \wedge e^3 \wedge e^4)= 2 e^1 \wedge e^2 \wedge e^3 \wedge e^4$ and thus we get the results of the Table for this case. The other cases can be handled in a similar way to complete the proof of the Table (\ref{coho}).
\end{proof}

\section{K\"ahler structures}

In this section we would like to search for (left) invariant K\"ahler structures, that is for symplectic forms $\omega$  on $\ggo$ such that there exists a complex structure $J$, which  is compatible with $\omega$, that  is 
$$\omega (JX, JY) = \omega (X,Y) \qquad \text{ for all } X, Y \in \ggo\qquad.$$
 To this end it is necessary to compare the results in Tables (\ref{complejas}) and (\ref{simplecticas}). We do not need to work with any complex structure but with a representative of its equivalence class. In fact, assume that there is a complex structure $J_1$ such that there is a symplectic structure $\omega$, which satisfies $\omega (J_1 X , J_1 Y) = \omega (X,Y)$ for all $X, Y \in \ggo$ and assume that  $J_1$ is equivalent to $J_2$. Thus there exists an automorphism $\sigma \in \Aut(\ggo)$ such that $J_2 = \sigma_{\ast}^{-^1} J_1 \sigma_{\ast}$. Then it holds
$$\omega (X , Y) = {\sigma^{\ast}}^{-^1}\sigma^{\ast} \omega (X , Y) = {\sigma^{\ast}}^{-^1}\omega (J_1 {\sigma^{\ast}}X,J_1\sigma^{\ast}Y)
=\omega (J_2 X,J_2Y).$$
Thus if there exists a symplectic structure $\omega$ compatible with $J_1$ then $\omega$ is also compatible with any $J_2$ which is equivalent to $J_1$, or equivalently if there is no symplectic structure compatible with $J_1$, then there is no symplectic structure compatible with $J_2$. 

Any such K\"ahler structure defines a pseudo-Riemannian metric $\phi$ on the Lie algebra in the following form: $\varphi(X,Y) = \omega(X, JY)$. The so defined  pseudo-Riemannian metric is compatible with the complex structure $J$. 

The following table shows the results of our investigation. For a given complex structure we describe the general form of a compatible symplectic structure when it exists.

\begin{prop} \label{kahler} Let $G$ be a solvable real Lie group of dimension four. Then the following Table shows for a given left invariant complex structure the general form of a compatible left invariant symplectic structure $\omega$ on $G$, when it exists (the so defined complex structure must be extended by imposing $J^2 = - \id$):

\begin{center} 
\small{
\begin{tabular}{|c|c|c|}\hline 
{\quad{\rm Case}\quad} & {{\rm \quad Complex structure \quad}} & {\rm K\"ahler structure}\\ \hline 
& & \\
{$\rr \rr_{3,0}$} &  {\quad $J e_1 = e_2,\, J e_3 = e_4$ \quad}& {$a_{12} e^1 \wedge e^2 + a_{34}e^3 \wedge e^4$} \\
& & {$a_{12}a_{34} \ne 0$}  \\ \hline
&& {$\quad a_{12}e^1 \wedge e^2+a_{13}e^1 \wedge e^3 + a_{14}e^1 \wedge e^4 - \quad $}\\
{$\rr \hh_{3}$} & {$J e_1 = e_2,\, J e_3 = e_4$}  & {$-a_{14}e^2 \wedge e^3+a_{13}e^2 \wedge e^4$}\\
 &&{$a_{13}^2+a_{14}^2\ne 0$} \\ \hline
 &&\\
{$\rr \rr_{3,0}$} & {$J e_1 = e_2,\, J e_3 = e_4$}& {$a_{12}e^1 \wedge e^2 + a_{34}e^3 \wedge e^4$} \\
&& {$a_{12}a_{34}\ne 0$} \\ \hline
 &&\\
{$\rr \rr'_{3,0}$} & {$J e_1 = e_4,\, J e_2 = e_3$}& {$a_{14}e^1 \wedge e^4 + a_{23}e^2 \wedge e^3$} \\
&& {$a_{14}a_{23}\ne 0$} \\ \hline
&& \\
{$\rr_2 \rr_2$} & {$J e_1 = e_2,\, J e_3 = e_4$}& {$a_{12}e^1 \wedge e^2 + a_{34}e^3 \wedge e^4$}\\
& & {$a_{12}a_{34}\ne 0$} \\ \hline 
&& \\
{$\rr'_2$} & {$J_1 e_1 = e_3,\, J_1 e_2 = e_4$} & {$a_{13-24}(e^1 \wedge e^3 - e^2 \wedge e_4)+\qquad \qquad \qquad \qquad$}\\
& &  {$+a_{14+23}(e^1 \wedge e^4 + e^2 \wedge e_3),\,a_{13-24}^2 +a_{14+23}^2 \ne 0$} \\\cline{2-3} & & \\
& {$J_2 e_1 = -e_2,\, J_1 e_3 = e_4$} & {$a_{12}( e^1 \wedge e^2) + a_{13-24}(e^1 \wedge e^3 - e^2 \wedge e_4)+$}\\
& &  {$+a_{14+23}(e^1 \wedge e^4 + e^2 \wedge e_3),\,a_{13-24}^2 +a_{14+23}^2 \ne 0$} \\\hline 
& & \\
{$\rr_{4,-1,-1}$} & {$Je_4 =  e_1,\, J e_2 = e_3$} &  {$a_{12}e^1 \wedge e^2 + a_{13} e^1 \wedge e^3+a_{14} e^1 \wedge e^4+$} \\
 & &{$-a_{13} e^2 \wedge e^4+ a_{12} e^3 \wedge e^4,\,a_{12}^2+a_{13}^2\ne 0$} \\ \hline
\end{tabular}} 
 \end{center}
 
\pagebreak

\begin{center}
\small{
\begin{tabular}{|c|c|c|}\hline
& & \\
{$\rr_{4,-1,1}$} & {$Je_4 =  e_2,\, J e_1 = e_3$} &  {$a_{12}e^1 \wedge e^2 + a_{13} e^1 \wedge e^3$} \\
 & &{$+a_{24} e^2 \wedge e^4- a_{12} e^3 \wedge e^4,\,a_{12}\ne 0$} \\ \hline
&&\\
{$\rr'_{4,0, \delta}$} & {$J_1 e_4 = e_1,\, J_1 e_2 = e_3,\, $}&
{$a_{14}e^1 \wedge e^4 + a_{23}e^2 \wedge e^3 $} \\
&{$J_2 e_4 = e_1,\, J_2 e_2 = -e_3$}& {$a_{14}a_{23}\ne 0$}\\   \hline 
& & \\
{$\dd_{4,1}$} & {$J_1 e_1 = e_4,\, J_1 e_2 = e_3$} &{ $a_{12-34}(e^1 \wedge e^2 - e^3 \wedge e^4)+e_{14} (e^1 \wedge e^4)$} \\
& & {$a_{12-34}\ne 0$} \\ \hline
&&\\
{$\dd_{4,2}$} & {$J_1e_4 = - e_2,\, J_1 e_1 = e_3$} & {$a_{14}e^1 \wedge e^4 + a_{14} e^2 \wedge e^3+ a_{24} e^2 \wedge e^4$}\\ & & {$a_{14}\ne 0$} \\ \cline{2-3}
&&\\
 &{$J_2e_4 =  -2e_1,\, J_2 e_2 = e_3$}&{$ a_{14} e^1 \wedge e^4+a_{23} e^2 \wedge e^3$}\\
 & &{$a_{14}a_{23}\ne 0$}\\ \hline
&&\\
{$\dd_{4,1/2}$} & {$J_1e_4 =  e_3,\, J_1 e_1 = e_2$} & {$a_{12-34} (e^1\wedge e^2-e^3 \wedge e^4)$}\\
& &{$a_{12-34}\ne 0$} \\\cline{2-3}
&&\\
 {}& {$J_2e_4 =  e_3,\, J_2 e_1 = -e_2$}& {$a_{12-34} (e^1\wedge e^2-e^3 \wedge e^4)$} \\ 
&&{$a_{12-34}\ne 0$}\\\hline
&{$J_1e_4 =  e_3,\, J_1 e_1 = e_2$}&\\
{$\dd'_{4,\delta}$} & {$J_2e_4 =  -e_3,\, J_2 e_1 = e_2$}& { $a_{12-\delta 34} (e^1\wedge e^2-\delta e^3 \wedge e^4)$}\\
& {$J_3e_4 =  -e_3,\, J_3 e_1 = -e_2$}& \\ 
{$\delta \ne 0$}&{$J_4e_4 =  e_3,\, J_4 e_1 = e_2$}&{$a_{12-34}\ne 0$} \\ \hline 
\end{tabular} }
\end{center}
\begin{center} 
{\rm Table} \ref{kahler} 
\end{center} 
\end{prop}
\begin{proof} We will give the proof in the case $\rr'_2$ since all cases should be handled in a simialr way. As we can see on the Table (\ref{complejas}) the complex structures on $\rr_2'$ are given by: $J_1 e_1 = e_3$, $J_1e_2 = e_4$;  and for the other type of complex structures, denoting $a_1\in \CC$ by $a_1 = \mu + i \nu$, with $\nu \ne 0$; we have   $J_{\mu,\nu} e_1 =\frac{\mu}{\nu} e_1 +(\frac{\nu^2 +\mu ^2}{\nu})e_2 $, $J_{\mu, \nu} e_3 = e_4$. On the other hand any sympletic structure has the form: $\omega = a_{12}(e^1 \wedge e^2) + a_{13-24}(e^1 \wedge e^3 - e^2 \wedge e_4)+a_{14+23}(e^1 \wedge e^4 + e^2 \wedge e_3)$, with $ a_{14+23}^2 + a_{13-24}^2 \ne 0$. Assuming that there exists a K\"ahler structure it  holds $\omega (JX,JY) = \omega (X,Y)$ for all $X,Y \in \ggo$ and this condition produces  equations on the coefficients of $\omega$ which should be verified in each case. 

So for $J_1$ we need to compute only the following:
$$\omega (e_1,e_2) = a_{12} = \omega (e_3, e_4)$$ 
and $$\omega(e_1,e_4) = a_{14+23} = \omega (e_3, -e_2)$$
Thus these equalities impose the condition $a_{12}=0$. And so any K\"ahler structure corresponding to $J_1$ has the form  $\omega = a_{13-24}(e^1 \wedge e^3 - e^2 \wedge e^4)+a_{14+23}(e^1 \wedge e^4 + e^2 \wedge e^3)$ with $a_{13-24}^2+a_{14+23}^2\ne 0$.

For the second case correspondig to $J_{\mu,\nu}$, by computing $\omega(e_2,e_4)$, $\omega(e_1,e_3)$, we get respectively:
$$
\begin{array}{lrcl}
{\text{i)}} & (1 + \frac1{\nu}) a_{24-13} & = & \frac{\mu}{\nu} a_{14+23} \\
&&&\\
{\text{ii)}} & (1 + \frac{\mu^2 + \nu^2}{\nu}) a_{24-13} & = & -\frac{\mu}{\nu} a_{14+23} \\
\end{array}
$$
By comparing i) and ii)  we get:
$$
\begin{array}{lrcl}
 &(1 + \frac1{\nu})  a_{24-13} & = & -(1 + \frac{\mu^2 + \nu^2}{\nu}) a_{24-13}  \\
\end{array}
$$
and  this equality  implies iii)$ a_{24-13}=0$ or iv)$1 + \frac1{\nu}+ 1 + \frac{\mu^2 + \nu^2}{\nu} =0$. As $a_{24-13}\ne 0$ (since in this case we would also get $a_{14+23}=0$ and this would be  a contradiction) then it must hold iv), that is  $1+ \mu^2 + \nu^2 =-2\nu$ and that implies $ \mu^2 =-2\nu- 1- \nu^2= -(\nu +1)^2$  and so that is possible only if $\mu =0$ and $\nu = -1$. For this complex structure $J$, given by $Je_1=-e_2$ $ J e_3 = e_4$,   it is not difficult to prove  that for any symplectic structure $\omega$ it always holds $\omega(JX,JY) $ = $ \omega(X,Y)$, that is, any symplectic structure on $\ggo$ is compatible with $J$. In this way we have completed the proof of the assertion.
\end{proof}

\end{document}